\documentclass[11pt,reqno]{amsart}
\usepackage{enumerate, latexsym, amsmath, amsfonts, amssymb, amsthm, color}
\def\pmod #1{\ ({\rm{mod}}\ #1)}
\def\Z{\Bbb Z}
\def\N{\Bbb N}

\def\l{\left}
\def\r{\right}
\def\bg{\bigg}
\def\({\bg(}
\def\){\bg)}
\def\t{\text}
\def\f{\frac}

\def\ls{\leqslant}
\def\gs{\geqslant}
\def\se {\subseteq}
\def\sm{\setminus}
\def\bi{\binom}

\def\eq{\equiv}

\def\da{\delta}

\def\Proof{\noindent{\it Proof}}

\theoremstyle{plain}
\newtheorem{theorem}{Theorem}

\newtheorem{lemma}{Lemma}
\newtheorem{corollary}{Corollary}
\newtheorem{conjecture}{Conjecture}
\theoremstyle{definition}

\theoremstyle{remark}
\newtheorem{remark}{Remark}

\allowdisplaybreaks

 \vspace{4mm}

\begin{document}

\hbox{Preprint, {\tt arXiv:2411.14308}}
\medskip

\title
[New results similar to Lagrange's four-square theorem]
{New results similar to Lagrange's four-square theorem}
\author
[Zhi-Wei Sun] {Zhi-Wei Sun}

\address{School of Mathematics, Nanjing
University, Nanjing 210093, People's Republic of China}
\email{zwsun@nju.edu.cn}

\keywords{Additive base, quadratic polynomial, polygonal number.
\newline \indent Supported
by the National Natural Science Foundation of China (grant 12371004).}

\begin{abstract} In this paper we establish some new results similar to Lagrange's four-square theorem.
For example, we prove that any integer $n>1$ can be written as
$w(5w+1)/2+x(5x+1)/2+y(5y+1)/2+z(5z+1)/2$ with $w,x,y,z\in\mathbb Z$.

Let $a$ and $b$ be integers with $a>0$, $b>-a$ and $\gcd(a,b)=1$. When $2\nmid ab$, we show that
any sufficiently large integer can be written as
$$\frac{w(aw+b)}2+\frac{x(ax+b)}2+\frac{y(ay+b)}2+\frac{z(az+b)}2$$
with $w,x,y,z$ nonnegative integers. When $2\mid a$ and $2\nmid b$,  we prove that any sufficiently large integer can be written as
$$w(aw+b)+x(ax+b)+y(ay+b)+z(az+b)$$
with $w,x,y,z$ nonnegative integers.
\end{abstract}
\maketitle

\section{Introduction}
\setcounter{lemma}{0}
\setcounter{theorem}{0}
\setcounter{corollary}{0}
\setcounter{remark}{0}
\setcounter{equation}{0}

For any integer $m\gs3$, the {\it polygonal numbers of order $m$} (or {\it $m$-gonal numbers}) are the nonnegative integers
$$p_m(n)=(m-2)\bi n2+n=\f{n((m-2)n-(m-4))}2\ \ (n=0,1,2,\ldots).$$
Those
$$\bar p_m(n)=p_m(-n)=\f{n((m-2)n+(m-4))}2\ \ (n=0,1,2,\ldots)$$
are called {\it second $m$-gonal numbers}, and those $p_m(x)$ with $x\in\Z$
are called {\it generalized $m$-gonal numbers}. Note that
\begin{gather*}p_3(x)=\f{x(x+1)}2,\ p_4(x)=x^2,\ p_5(x)=\f{x(3x-1)}2,
\\ p_6(x)=x(2x-1)=p_3(2x-1),\ p_7(x)=\f{x(5x-3)}2,\ p_8(x)=x(3x-2).
\end{gather*}

In 1638 Fermat claimed that each $n\in\N=\{0,1,2,\ldots\}$ is a sum of $m$ $m$-gonal numbers;
this was confirmed by
Lagrange, Gauss and Cauchy for the cases $m=4$, $m=3$, and $m\gs5$, respectively (cf. \cite{N87,N96}).
Lagrange's four-square theorem states that any $n\in\N$ can be written as $w^2+x^2+y^2+z^2$
with $w,x,y,z\in\N$.
In 1830 Legendre refined Cauchy's polygonal number theorem by showing that for any integer $m\gs 5$
we can write each integer $n\gs 28(m-2)^3$ as $$p_{m}(w)+p_{m}(x)+p_{m}(y)+p_{m}(z)+\da_m(n),$$
where $w,x,y,z\in\N$, $\da_m(n)=0$ if $2\nmid m$, and $\da_m(n)\in\{0,1\}$ if $2\mid m$.
M. B. Nathanson (\cite{N87} and \cite[p.\,33]{N96}) simplified the proofs of Cauchy's and Legendre's theorems.
As a supplement to Legendre's theorem, X.-Z. Meng and the author \cite{MS} prove that
 for any integer $m\gs 6$ with $m\eq 2\pmod4$ each integer
$n\gs28(m-2)^3$ can be expressed as
$p_{m}(w)+p_{m}(x)+p_{m}(y)+p_{m}(z)$ with $w,x,y,z\in\N$, and that for any integer $m\gs8$ with $4\mid m$
there are infinitely many positive integers not of the form
$p_{m}(w)+p_{m}(x)+p_{m}(y)+p_{m}(z)$ with $w,x,y,z\in\N$.

In additive number theory, a subset $A$ of $\N$ is called a {\it base of order $n$} if
$$nA=\{a_1+\cdots+a_n:\ a_1,\ldots,a_n\in A\}$$
coincides with $\N$. When $\N\sm nA$ is finite, we call $A$ an {\it asymptotic base of order $n$}.
In this terminology, for each integer $m\gs3$ the set $\{p_m(n):\ n\in\N\}$
is a base of order $m$, and $\{p_m(n):\ n\in\N\}$ is an asymptotic base of order four if and only if
$m=4$ or $4\nmid m$.

In 1994 R. K. Guy \cite{G} observed that each $n\in\N$ is a sum of three generalized pentagonal numbers
(see also the comments around \cite[(1.4)]{S15}), that is, the set $\{p_5(x):\ x\in\Z\}$ is a base of order three.
In 2016 the author \cite{S16} proved that $\{p_8(x):\ x\in\Z\}$
is a base of order four; moreover, each $n\in\Z^+=\{1,2,3,\ldots\}$
can be written a sum of four generalized octagonal numbers one of which is odd.
This is quite similar to Lagrange's four-square theorem.
For any integer $m\gs5$, J. Li and H. Wang \cite{LW} used Jacobi forms to obtain a formula for
$$R_m(n)=|\{(w,x,y,z)\in\Z^4:\ p_m(w)+p_m(x)+p_m(y)+p_m(z)=n\}|$$
in terms of Hurwitz class numbers. Consequently,
for any integer $n\gs (m-2)^3/8+m/2+1$,
we have $R_m(n)>0$ if $m$ or $n$ is odd.

For any $a\in\Z^+$ and $b\in\Z$ with $b>-a$ and $a\eq b\pmod2$, we clearly have
$$\l\{\f{x(ax+b)}2:\ x\in\N\r\}\se\N.$$
In \cite{S18,S20} the author investigated those tuples $(a,b,c,d,e,f)\in\Z^6$
with $a,c,e\in\Z^+$ and $a+b,c+d,e+f\in 2\Z^+$
such that each $n\in\N$ can be written as
$$\f{x(ax+b)}2+\f{y(cy+d)}2+\f{z(ez+f)}2$$
with $x,y,z\in\N$. For example, the author \cite{S20} conjectured that
each $n\in\N$ can be written as $x(x+1)/2+y(3y+1)/2+z(5z+1)/2$ with $x,y,z\in\N$, and proved that
$$\l\{\f{x(x+1)}2+\f{y(3y+1)}2+\f{z(5z+1)}2:\ x,y,z\in\Z\r\}=\N.$$

Motivated by the above work, we establish some new results
similar to Lagrange's four-square theorem.

\begin{theorem} \label{Th1.1} Let $a\in\Z^+$ and $b\in\Z$ with $b>-a$, $2\nmid ab$ and $\gcd(a,b)=1$.
Then the set $\{x(ax+b)/2:\ x\in\N\}$ is an asymptotic base of order four. Moreover, for any
$n\in\N$ with
\begin{equation}\label{n>}
n>2a(7a(a-1)+1)+(7a-2)b+(2a(2a-1)+b)\sqrt{6(2a(a-1)+b)},
\end{equation}
there are $w,x,y,z\in\N$ such that
$$n=\f{w(aw+b)}2+\f{x(ax+b)}2+\f{y(ay+b)}2+\f{z(az+b)}2.$$
\end{theorem}

Via Theorem \ref{Th1.1} and some easy numerical computations, we obtain the following corollary
involving second pentagonal numbers and second heptagonal numbers.

\begin{corollary}\label{Cor-p57}  We have
\begin{align*}& \l\{\bar p_5(w)+\bar p_5(x)+\bar p_5(y)+\bar p_5(z):\ w,x,y,z\in\N\r\}
\\= &\ \N\sm\{1,3,5,10,12,20,25,27,38,53,65,153,165\}.
\end{align*}
Also,
any integer $n\gs 786$ can be written as
$\bar p_7(w)+\bar p_7(x)+\bar p_7(y)+\bar p_7(z)$
with $w,x,y,z\in\N$.
\end{corollary}
\begin{remark} \label{Rem-p57} By Legendre's refinement of Cauchy's polygonal number theorem,
any integer $n\gs 28\times (5-2)^3=756$ is a sum of four pentagonal numbers,
and integer $n\gs 28\times (7-2)^3=3500$ is a sum of four heptagonal numbers.
Combining this with direct checks for small numbers, we find that
\begin{align*}\l\{p_5(w)+p_5(x)+p_5(y)+p_5(z):\ w,x,y,z\in\N\r\}= \N\sm\{9,21,31,43,55,89\}
\end{align*}
and
$$\{p_7(w)+p_7(x)+p_7(y)+p_7(z):\ w,x,y,z\in\N\}\supseteq\{n\in\N:\ n\gs 1004\}.$$
\end{remark}

Theorem \ref{Th1.1} with $(a,b)=(5,\pm1)$ has the following consequence.

\begin{corollary} \label{Cor5-1} Each integer $n>1$ can be written as
$$\f{w(5w+1)}2+\f{x(5x+1)}2+\f{y(5y+1)}2+\f{z(5z+1)}2$$
with $w,x,y,z\in\Z$. Moreover, any integer $n\gs 776$ can be written as
$$\f{w(5w+1)}2+\f{x(5x+1)}2+\f{y(5y+1)}2+\f{z(5z+1)}2$$
with $w,x,y,z\in\N$, and
any integer $n\gs 884$ can be written as
$$\f{w(5w-1)}2+\f{x(5x-1)}2+\f{y(5y-1)}2+\f{z(5z-1)}2$$
with $w,x,y,z\in\N$.
\end{corollary}
\begin{remark} Via a computer we can easily verify that the set
 $$\l\{\f{w(5w+1)}2+\f{x(5x+1)}2+\f{y(5y+1)}2+\f{z(5z+1)}2:\ w,x,y,z\in\Z\r\}$$
 contains $2,3,\ldots,775.$
So the first assertion in Corollary \ref{Cor5-1} follows from the second assertion in the corollary.
\end{remark}

For convenience, for $a\in\Z^+$ and $b\in\Z$ with  $a\eq b\pmod 2$, we set
$$S_{a,b}=\l\{\f{w(aw+b)}2+\f{x(ax+b)}2+\f{y(ay+b)}2+\f{z(az+b)}2:\ w,x,y,z\in\Z\r\}.$$

Via Theorem \ref{Th1.1} and some easy computations via a computer, we have  the following corollary.

\begin{corollary}\label{Cor1.1} We have
$$S_{7,1}=\N\sm\{1,2,5\},\ S_{7,3}=\N\sm\{1,3,25\},
\ S_{7,5}=\N\sm\{5,23\}.$$
Also,
\begin{align*}S_{9,1}&=\N\sm\{1,2,3,6,7,11,35,37\},
\\ S_{9,5}&=\N\sm\{1,3,5,10,12,31,67\},
\\S_{9,7}&=\N\sm\{5,6,7,15,29,65\}.
\end{align*}
\end{corollary}

Now we state our second theorem.

\begin{theorem}\label{Th1.2} Let $a\in\Z^+$ and $b\in\Z$ with $b>-a$, $2\mid a$ and $\gcd(a,b)=1$.
Then the set $\{x(ax+b):\ x\in\N\}$ is an asymptotic base of order four. Moreover, if $n$
is an odd integer greater than
\begin{equation*}
4a(7a(a-1)+1)+2(7a-2)b+2(2a(2a-1)+b)\sqrt{6(2a(a-1)+b)},
\end{equation*}
 or $n$ is an even integer with $n/4$ greater than
$$28a^3-14a^2+a+(7a-1)b+\l(2a(4a-1)+b\r)\sqrt{3(2a(2a-1)+b)},$$
 then we can write $n$ as
$$w(aw+b)+x(ax+b)+y(ay+b)+z(az+b)$$
with $w,x,y,z\in\N$.
\end{theorem}

For $a\in\Z^+$ and $b\in\Z$, we define
$$T_{a,b}=S_{2a,2b}=\{w(aw+b)+x(ax+b)+y(ay+b)+z(az+b):\ w,x,y,z\in\Z\}.$$

\begin{corollary}\label{Cor1.2} We have
$$ T_{4,1}=\N\sm\{1,2,4,7,30\},\ \max(\N\sm T_{6,1})=168,\ \max(\N\sm T_{6,5})=182.$$
Moreover, any integer $n\gs 2283$ can be written as
$w(4w+1)+x(4x+1)+y(4y+1)+z(4z+1)$ with $w,x,y,z\in\N$, and each integer $n>4500$ can be written as
$w(4w-1)+x(4x-1)+y(4y-1)+z(4z-1)$ with $w,x,y,z\in\N$.
\end{corollary}
\begin{remark} Meng and Sun \cite[Corollary 1.1]{MS}
showed that $T_{4,3}=\{p_{10}(w)+p_{10}(x)+p_{10}(y)+p_{10}(z):\ w,x,y,z\in\Z\}$
coincides with $\N\sm\{5,6,26\}$.
\end{remark}

Now we present our third theorem.

\begin{theorem} \label{Th1.3} Let $a\in\Z^+$ and $b\in\Z$ with $b>-a$, $2\mid b$ and $\gcd(a,b)=1$. Then, for any integer $n\not\eq0\pmod4$ greater than
$$4a(7a(a-1)+1)+2(7a-2)b+2(2a(2a-1)+b)\sqrt{6(2a(a-1)+b)},$$
there are $w,x,y,z\in\N$ such that
$$n=w(aw+b)+x(ax+b)+y(ay+b)+z(az+b).$$
\end{theorem}
\begin{remark} Let $a\in\Z^+$ and $b\in\Z$ with $b>-a$, $2\mid b$ and $\gcd(a,b)=1$.
As $a+b$ is odd, for each $n\in\N$ either $n$ or $n-(a+b)$ is not divisible by $4$.
So Theorem \ref{Th1.3} implies that the set $\{x(ax+b):\ x\in\N\}$ is an asymptotic base of order five.
\end{remark}

\begin{corollary} We have
$$\{n\in\N:\ 4\nmid n\ \t{and}\ n\not\in T_{5,2}\}=\{1,2,5,11,18\}$$
and
$$\{n\in\N:\ 4\nmid n\ \t{and}\ n\not\in T_{5,4}\}=\{5,6,7,17\}.$$
\end{corollary}
\begin{remark} Sun \cite{S16} proved that $T_{3,2}=\N$. It seems that
both $\N\sm T_{5,2}$ and $\N\sm T_{5,4}$ contain infinitely many multiples of $4$.
\end{remark}

In 2020 D. Krachun and the author \cite{KS} proved that
$$\{2p_5(w)+p_5(x)+p_5(y)+p_5(z):\ w,x,y,z\in\N\}=\N,$$
which was previously conjectured by the author \cite{S16}.
Motivated by this, we establish the following theorem.

\begin{theorem} \label{Th2} Let $a\in\Z^+$ and $b\in\Z$ with $b>-a$ and $\gcd(a,5b)=1$.
Then, for any integer
\begin{equation}\label{25}n\gs 550a(60a^2+b)+25(120a^2+b)\sqrt{2(60a^2+b)}.
\end{equation}
with $2\mid n$ if $2\nmid ab$, there are $w,x,y,z\in\N$ such that
\begin{equation} 2w(aw+b)+x(ax+b)+y(ay+b)+z(az+b)=n.
\end{equation}
\end{theorem}
\begin{remark} We guess that $\gcd(a,5b)$ in Theorem \ref{Th2} can be replaced by $\gcd(a,b)$,
but we are unable to prove this.
\end{remark}

Theorem \ref{Th2} in the case $2\nmid ab$ yields the following result.

\begin{corollary}\label{Cor2} Let $a\in\Z^+$ and $b\in\Z$ with $b>-a$, $2\nmid ab$ and $\gcd(a,5b)=1$.
Then any integer
\begin{equation*}\gs 275a(60a^2+b)+\f{25}2(120a^2+b)\sqrt{2(60a^2+b)}.
\end{equation*}
can be written as
\begin{equation} w(aw+b)+\f{x(ax+b)}2+\f{y(ay+b)}2+\f{z(az+b)}2
\end{equation}
with $w,x,y,z\in\N$.
\end{corollary}

In the case $(a,b)=(3,1)$, Corollary \ref{Cor2} yields the following concrete result.

\begin{corollary} We have
\begin{align*}&\l\{2\bar p_5(w)+\bar p_5(x)+\bar p_5(y)+\bar p_5(z):\, w,x,y,z\in\N\r\}
\\&\qquad =\N\sm\{1,3,5,12,27,50\},
\end{align*}
\end{corollary}

Theorem \ref{Th2} with $(a,b)=(2,\pm1)$ yields the following corollary.

\begin{corollary} We have
$$\{2w(2w-1)+x(2x-1)+y(2y-1)+z(2z-1):\ w,x,y,z\in\N\}=\N\sm\{11,26,50\}.$$
Also, any integer $n\gs90$ can be written as $2w(2w+1)+x(2x+1)+y(2y+1)+z(2z+1)$
with $w,x,y,z\in\N$.
\end{corollary}

Meng and Sun \cite[Theorem 1.3]{MS} proved that
for each integer $m\gs5$, any integer $n\gs 924(m-2)^3$ can be written as
$3p_m(w)+p_m(x)+p_m(y)+p_m(z)$ with $w,x,y,z\in\N$. We extend this as follows.

\begin{theorem}\label{Th3} Let $a\in\Z^+$ and $b\in\Z$ with $b>-a$ and $\gcd(a,b)=1$.
Let $n\in\N$ with $ab$ or $n$ even. Suppose that
\begin{equation}\label{n6}\f n6>66a^2(9a-1)+a+(33a-1)b+(2a(18a-1)+b)\sqrt{15(2a(9a-1)+b)}.
\end{equation}
Then there are $w,x,y,z\in\N$ such that
$$3w(aw+b)+x(ax+b)+y(ay+b)+z(az+b)=n.$$
\end{theorem}

Theorem \ref{Th3} in the case $2\nmid ab$ yields the following result.

\begin{corollary} \label{Cor1.9} Let $a$ and $b$ be odd integers with $a>0$,  $b>-a$ and $\gcd(a,b)=1$.
Then, any integer $n$ with
\begin{equation}\label{n3}\f n3>66a^2(9a-1)+(33a-1)b+(2a(18a-1)+b)\sqrt{15(2a(9a-1)+b)}.
\end{equation}
can be written as
$$3\f{w(aw+b)}2+\f{x(ax+b)}2+\f{y(ay+b)}2+\f{z(az+b)}2$$
with $w,x,y,z\in\N$.
\end{corollary}

Applying Corollary \ref{Cor1.9} with some concrete odd integers $a$ and $b$, we obtain the following corollary.

\begin{corollary} {\rm (i)} We have
$$\l\{3\bar p_5(w)+\bar p_5(x)+\bar p_5(y)+\bar p_5(z):\ w,x,y,z\in\N\r\}=\N\sm\{1,3,5,18,31\}.$$

{\rm (ii)} Any integer $n>1$ can be written as
$$3\f{w(5w+1)}2+\f{x(5x+1)}2+\f{y(5y+1)}2+\f{z(5z+1)}2$$
with $w,x,y,z\in\Z$. Moreover, each integer $n>151$ can be written as
$$3\f{w(5w-1)}2+\f{x(5x-1)}2+\f{y(5y-1)}2+\f{z(5z-1)}2$$
with $w,x,y,z\in\N$, and any integer $n>682$ can be written as
$$3\f{w(5w+1)}2+\f{x(5x+1)}2+\f{y(5y+1)}2+\f{z(5z+1)}2$$
with $w,x,y,z\in\N$.

{\rm (iii)} $1,2,5,14$ are the only natural numbers not of the form
$$3\f{w(7w+1)}2+\f{x(7x+1)}2+\f{y(7y+1)}2+\f{z(7z+1)}2\ \ (w,x,y,z\in\Z),$$
and $1,3,14$ are the only natural numbers not of the form
$$3\f{w(7w+3)}2+\f{x(7x+3)}2+\f{y(7y+3)}2+\f{z(7z+3)}2\ \ (w,x,y,z\in\Z).$$
\end{corollary}
\begin{remark} Meng and Sun \cite[Corollary 1.3]{MS} proved that
\begin{align*}&\ \{3p_7(w)+p_7(x)+p_7(y)+p_7(z):w,x,y,z\in\N\}
\\=&\ \N\sm\{13,16,27,31,33,49,50,67,87,178,181,259\}
\end{align*}
and
$$\{3p_9(w)+p_9(x)+p_9(y)+p_9(z):w,x,y,z\in\Z\}=\N\sm\{17\}.$$
\end{remark}

Theorem \ref{Th3} with $(a,b)=(4,\pm1)$ yields the following corollary.
\begin{corollary}\label{Cor41} We have
\begin{align*}&\{3w(4w+1)+x(4x+1)+y(4y+1)+z(4z+1):\ w,x,y,z\in\Z\}
\\&\qquad=\N\sm\{1,2,4,7,16\}.
\end{align*}
Moreover, any integer $n>1204$ can be written as
$3w(4w-1)+x(4x-1)+y(4y-1)+z(4z-1)$ with $w,x,y,z\in\N$, and
any integer $n>2104$ can be written as
$3w(4w+1)+x(4x+1)+y(4y+1)+z(4z+1)$ with $w,x,y,z\in\N$.
\end{corollary}

\begin{remark} In contrast with Corollary \ref{Cor41}, Meng and Sun \cite[Corollary 1.3]{MS} proved that
$$ \{3p_{10}(w)+p_{10}(x)+p_{10}(y)+p_{10}(z):w,x,y,z\in\Z\}=\N\sm\{16,19\}.$$
\end{remark}

We are going to prove Theorems \ref{Th1.1}-\ref{Th1.3},
\ref{Th2}, \ref{Th3} in Sections 2, 3, 4 respectively.
In Section 5, we give some additional remarks and pose some open conjectures.

\section{Proofs of Theorems \ref{Th1.1} and \ref{Th1.3}}
\setcounter{lemma}{0}
\setcounter{theorem}{0}
\setcounter{corollary}{0}
\setcounter{remark}{0}
\setcounter{conjecture}{0}
\setcounter{equation}{0}

The following lemma is essentially \cite[Lemma 2.4]{MS}.

\begin{lemma} \label{Lem2.1} Let $c,d\in\N$ with $c\eq d\pmod 2$ and $4\nmid c$. If
\begin{equation}\label{cd}4c>d^2\ \t{and}\ 3c<d^2+2d+4,
\end{equation}
then there are $w,x,y,z\in\N$ such that
\begin{equation}\label{2.1} w^2+x^2+y^2+z^2=c\ \t{and}\ w+x+y+z=d.
\end{equation}
\end{lemma}
\begin{remark} Lemma \ref{Lem2.1} in the case $2\nmid cd$ is Cauchy's Lemma
(cf. \cite[p.\,31]{N96}) used to prove Cauchy's polygonal number theorem.
\end{remark}

\begin{lemma}\label{Lem2.2} Suppose that $ac+bd=n\gs2a$ with $a\in\Z^+$, $b,c,d\in\Z$ and $b>-a$.
Then
\begin{equation}\label{<c<}d>2\ \ \t{and}\ \ \f{d^2}4<c<\f{d^2+2d+4}3
\end{equation}
if and only if
\begin{equation}\label{<d<}\f{\sqrt{3a(n-a+b)+\f 94b^2}-a-3b/2}a<d<\f 2a(\sqrt{an+b}-b).
\end{equation}
\end{lemma}
\Proof. Clearly,
$$3a(n-a+b)+\f94b^2\gs 3a^2+3ab+\f 94b^2=2a^2+\l(a+\f{3b}2\r)^2\gs2.$$
Thus $d>2$ if \eqref{<d<} holds.

Now assume $d>2$. Then $ad+2b>a(d-2)>0$. Note that
$$n>a\f{d^2}4+bd\iff 4an>(ad+2b)^2-4b^2\iff d<\f2a\l(\sqrt{an+b^2}-b\r).$$
As
$$ad+a+\f{3b}2>\f 32(a+b)>0,$$
we have
\begin{align*}&\ n<\f a3(d^2+2d+4)+bd
\\\iff&\ 3an<a^2d^2+ad(2a+3b)+4a^2=\l(ad+a+\f{3b}2\r)^2+3a(a-b)-\f{9}4b^2
\\\iff&\sqrt{3a(n-a+b)+\f94b^2}<ad+a+\f{3b}2
\\\iff&d>\f{\sqrt{3a(n-a+b)+9b^2/4}-a-3b/2}a.
\end{align*}
As $n=ac+bd$, by the above, $d^2/4<c<(d^2+2d+4)/3$ if and only if \eqref{<d<} holds.

In view of the above, we have completed the proof. \qed

Clearly, Theorems \ref{Th1.1} and \ref{Th1.3} follow from the following theorem.

\begin{theorem} \label{Th2.1} Let $a\in\Z^+$ and $b\in\Z$ with $b>-a$ and $\gcd(a,b)=1$. Let $n\in\N$ with
\begin{equation}\label{n>}
n>4a(7a(a-1)+1)+2(7a-2)b+2(2a(2a-1)+b)\sqrt{6(2a(a-1)+b)}.
\end{equation}
Suppose that $2\nmid ab$ and $2\mid n$, or $2\mid a$ and $2\nmid n$, or $2\mid b$ and $4\nmid n$.
Then we can write $n$ as
$$w(aw+b)+x(ax+b)+y(ay+b)+z(az+b)$$
with $w,x,y,z\in\N$.
\end{theorem}
\Proof. Note that
$$2a(2a-1)+b>2a(a-1)+b\gs 2a(a-1)-a+1=a(2a-3)+1\gs0.$$
If $a>1$, then by \eqref{n>} we have
$$\f n2>2(7a^2-7a+1)+(7a-2)b>2(7a^2-7a+1)-a(7a-2)=a(7a-12)+2>a.$$
When $a=1$, we have $b\gs 1$ and hence
$$\f n2>2(7a^2-7a+1)+(7a-2)b\gs 14a^2-14a+2+7a-2=14a^2-7a>a.$$
Thus we always have $n>2a$.

In view of Lemma \ref{Lem2.1}, we want to find positive integers $c$ and $d$ with $c\eq d\pmod2$
and $4\nmid c$ such that
$$\f{d^2}4<c<\f{d^2+2d+4}3\ \t{and}\ n=ac+bd.$$

Motivated by Lemma \ref{Lem2.2}, we consider the interval
\begin{equation}\label{I}I=\l(\f{\sqrt{3a(n-a+b)+9b^2/4}-a-3b/2}a,\ \f2a\l(\sqrt{an+b^2}-b\r)\r).
\end{equation}
Clearly the length of $I$ is greater than $2a$ if and only if
\begin{equation}\label{an}2\sqrt{an+b^2}-a(2a-1)-\f b2>\sqrt{3a(n-a+b)+\f 94b^2}.
\end{equation}
As $2a(2a-1)+b>0$, we have
\begin{align*}&\ 2\sqrt{an+b^2}\gs a(2a-1)+\f b2
\\\iff&\ 4(an+b^2)\gs a^2(2a-1)^2+(2a-1)ab+\f{b^2}4.
\end{align*}
If $b\gs a(2a-1)$, then
$$4b^2\gs a^2(2a-1)^2+a(2a-1)b+2b^2>a^2(2a-1)^2+(2a-1)ab+\f{b^2}4.$$
When $b<a(2a-1)$, we have $a>1$ and
\begin{align*}\f n2&>2a(7a^2-7a+1)+(7a-2)b
\\&> 2a(7a^2-7a+1)-a(7a-2)=a(14a^2-21a+4)
\\&>a\f{(2a-1)^2}4
\end{align*}
and hence
$$4an>2a^2(2a-1)^2>a^2(2a-1)^2+a(2a-1)b.$$
Thus the left-hand side of \eqref{an} is nonnegative. Let $t=\sqrt{an+b^2}$. Then
\begin{align*}&\ \eqref{an}\ \t{holds}
\\\iff&\ \l(2t-a(2a-1)-\f b2\r)^2>3a(n-a+b)+\f 94b^2
\\\iff&\ 4t^2-2t(2a(2a-1)+b)+\l(a(2a-1)+\f b2\r)^2>3a(n-a+b)+\f 94b^2
\\\iff&\ (t-(2a(2a-1)+b))^2>6a^2(2a(a-1)+b),
\end{align*}
and
\begin{align*}&\ \l(2a(2a-1)+b+\sqrt{6a^2(2a(a-1)+b)}\r)^2
\\&\ -2(2a(2a-1)+b)a\sqrt{6(2a(a-1)+b)}
\\=&\ (2a(2a-1)+b)^2+6a^2(2a(a-1)+b)
\\=&\ 4a^2(7a^2-7a+1)+14a^2b-4ab+b^2
\\<&\ b^2+a(n-2(2a(2a-1)+b)\sqrt{6(2a(a-1)+b)})
\end{align*}
by \eqref{n>}. Thus
$$t=\sqrt{an+b^2}>2a(2a-1)+b+\sqrt{6a^2(2a(a-1)+b)}$$
and hence \eqref{an} does hold.

We claim that there are $c,d\in\Z$ with $c\eq d\pmod 2$ and $4\nmid c$ such that $d\in I$ and $ac+bd=n$.
To prove this, we distinguish two cases.

{\it Case} 1. $2\nmid ab$ and $2\mid n$.

In this case, $\gcd(b,2a)=1$. As the length of the interval $I$ is greater than $2a$, there is an integer $d\in I$ with $n-bd\eq a\pmod{2a}$. As $2\nmid a$ and $2\mid n$,
we must have $2\nmid d$. Write $n-bd =ac$ with $c$ an odd integer.

{\it Case} 2. $2\mid a$ and $2\nmid n$, or $2\mid b$ and $4\nmid n$.

In this case, $a-b$ is odd since $\gcd(a,b)=1$. Since $\gcd(a-b,2a)=1$ and the length of the interval $I$ is greater than $2a$, there is an integer $d\in I$ with $(b-a)d\eq n\pmod{2a}$.
Write $n=(b-a)d+2ac_0$ with $c_0\in\Z$. Then $n=ac+bd$ with $c=2c_0-d\eq d\eq n\pmod2$.
If $2\nmid n$, then we have $2\nmid cd$. If $2\mid b$ and $n\eq2\pmod 4$,
then $ac\eq ac+bd=n\eq2\pmod4$ and hence $c\eq2\pmod 4$.

Since $d\in I$ in either case, we have \eqref{<c<} by Lemma \ref{Lem2.2}. In light of
Lemma \ref{Lem2.1}, there are $w,x,y,z\in\N$ satisfying \eqref{2.1}. Therefore
\begin{align*}n&=ac+bd=a(w^2+x^2+y^2+z^2)+b(w+x+y+z)
\\&=w(aw+b)+x(ax+b)+y(ay+b)+z(az+b).
\end{align*}
This concludes our proof. \qed

\begin{theorem} \label{Th2.2} Let $a\in\Z^+$ and $b\in\Z$ with $2\mid a$, $\gcd(a,b)=1$ and $b>-a$. Let $n\in\N$ with $2\mid n$ and
\begin{equation}\label{n/4}\begin{aligned}
\f n4>&\ 28a^3-14a^2+a+(7a-1)b+\l(2a(4a-1)+b\r)\sqrt{3(2a(2a-1)+b)}.
\end{aligned}\end{equation}
Then there are $w,x,y,z\in\N$ such that
$$n=w(aw+b)+x(ax+b)+y(ay+b)+z(az+b).$$
\end{theorem}
\Proof. As $b>-a$, we have
$$\f n4>28a^3-14a^2+a-a(7a-1)=8a^3+a(20a^2-21a+2)>8a^3.$$

For the interval $I$ given by \eqref{I}, its length is greater than $4a$
 if and only if
\begin{equation}\label{4a}2\sqrt{an+b^2}-a(4a-1)-\f b2>\sqrt{3a(n-a+b)+\f 94b^2}.
\end{equation}
As $2a(4a-1)+b>2a(4a-1)-a>0$, we have
\begin{align*}&\ 2\sqrt{an+b^2}\gs a(4a-1)+\f b2
\\\iff&\ 4(an+b^2)\gs a^2(4a-1)^2+(4a-1)ab+\f{b^2}4.
\end{align*}
If $b\gs a(4a-1)$, then
$$4b^2\gs a^2(4a-1)^2+a(4a-1)b+2b^2>a^2(4a-1)^2+(4a-1)ab+\f{b^2}4.$$
When $b<a(4a-1)$, we have
$$4an>4a\times 8a^3>2a^2(4a-1)^2>a^2(4a-1)^2+a(4a-1)b.$$
Thus the left-hand side of \eqref{4a} is nonnegative. Let $t=\sqrt{an+b^2}$. Then
\begin{align*}&\ \eqref{4a}\ \t{holds}
\\\iff&\ \l(2t-a(4a-1)-\f b2\r)^2>3a(n-a+b)+\f 94b^2
\\\iff&\ 4t^2-2t(2a(4a-1)+b)+\l(a(4a-1)+\f b2\r)^2>3a(n-a+b)+\f 94b^2
\\\iff&\ (t-(2a(4a-1)+b))^2>12a^2(2a(2a-1)+b),
\end{align*}
and
\begin{align*}&\ \l(2a(4a-1)+b+\sqrt{12a^2(2a(2a-1)+b)}\r)^2
\\&\ -4a(2a(4a-1)+b)\sqrt{3(2a(2a-1)+b)}
\\=&\ (2a(4a-1)+b)^2+12a^2(2a(2a-1)+b)
\\=&\ b^2+28a^2(4a^2-2a+b)+4a^2-4ab
\\<&\ b^2+an-4a(2a(4a-1)+b)\sqrt{3(2a(2a-1)+b)}
\end{align*}
by \eqref{n/4}. Thus
$$t=\sqrt{an+b^2}>2a(4a-1)+b+\sqrt{12a^2(2a(2a-1)+b)}$$
and hence \eqref{4a} holds. So the length of the interval $I$ is greater than $4a$.

As $\gcd(b,4a)=1$, there is an integer $d\in I$ with $n-bd\eq 2a\pmod{4a}$. Write $n-bd =ac$ with
$c\eq2\pmod 4$. As both $a$ and $n$ are even but $b$ is odd, we must have $2\mid d$.
As $d\in I$, we have \eqref{<c<} by Lemma \ref{Lem2.2}. Now, by Lemma \ref{Lem2.1}
there are $w,x,y,z\in\N$ satisfying \eqref{2.1}.
Hence
$$n=ac+bd=w(aw+b)+x(ax+b)+y(ay+b)+z(az+b).$$
This concludes our proof. \qed

\medskip
\noindent {\tt Proof of Theorem \ref{Th1.2}}. As $2\mid a$ and $2\nmid b$,
combining Theorems \ref{Th2.1} and \ref{Th2.2} we immediately obtain the desired result. \qed

\section{Proof of Theorem \ref{Th2}}
\setcounter{lemma}{0}
\setcounter{theorem}{0}
\setcounter{corollary}{0}
\setcounter{remark}{0}
\setcounter{equation}{0}
\setcounter{conjecture}{0}

We need the following lemma which is \cite[Lemma 2.1]{KS}.

\begin{lemma}\label{Lem4.1} Any even $n\in\N$ not in the set
\begin{equation}\label{E}E=\{5^{2k+1}m:\ k,m\in\N \ \t{and}\ m\eq\pm2\pmod5\}
\end{equation}
can be written as $x^2+y^2+z^2+(x+y+z)^2/2$ with $x,y,z\in\Z$.
\end{lemma}

Our proof of Theorem \ref{Th2} extends the method in \cite{KS}.

\medskip
\noindent{\tt Proof of Theorem \ref{Th2}}. Note that $\gcd(a+b,a)=\gcd(b,a)=1$
and $\gcd(5(a+b),2a)=\gcd(a+b,2)$ divides $n$. So the congruence $5(a+b)x\eq n\pmod{2a}$
is solvable (over $\Z$). When $n=5q$ with $q\in\Z$, we choose $\da\in\{1,-1\}$ such that
$4aq+b^2+\da\not\eq\pm2\pmod5$, and hence the congruence $(2ax+b)^2\eq 4aq+b^2+\da\pmod 5$
is solvable. Thus, by the Chinese Remainder Theorem, among any $10a$ consecutive integers, there is an integer $B$ such that
$5(a+b)B\eq n\pmod{2a}$, and
$$a(q-bB-aB^2)\eq 4a(aB^2+bB-q)=(2aB+b)^2-(4aq+b^2)\eq \da\pmod 5$$ if $n=5q$ with $q\in\N$.
Note that $(n-5bB)/a-5B^2$ is an even integer. If $n=5q$ with $q\in\N$, then
$$\f15\l(\f{n-5bB}a-5B^2\r)=\f{a(q-bB-aB^2)}{a^2}\eq\f{\da}{a^2}\not\eq0,\pm2\pmod5.$$
So $(n-5bB)/a-5B^2$ is an even number not in the set $E$ given by \eqref{E}.

Suppose $B\in\Z^+$. Then $2aB+b\gs a+b>0$, $12aB+5b\gs5(a+b)>0$, and
\begin{align*}&0\ls \f{n-5bB}a-5B^2\ls B^2
\\\iff&5bB+5aB^2\ls n\ls 5bB+6aB^2
\\\iff&4an\gs 5(2aB+b)^2-5b^2\ \t{and}\ 24an+25b^2\ls (12aB+5b)^2
\\\iff&\sqrt {\f45an+b^2}\gs 2aB+b\ \t{and}\ \sqrt {24an+25b^2}\ls 12aB+5b
\\\iff&\f{\sqrt{24an+25b^2}-5b}{12a}\ls B\ls \f{\sqrt{4an/5+b^2}-b}{2a}.
\end{align*}

Observe that
\begin{align*}&\f{\sqrt{4an/5+b^2}-b}{2a}-\f{\sqrt{24an+25b^2}-5b}{12a}\gs 10a
\\\iff&6\sqrt{\f 45an+b^2}\gs 120a^2+b+\sqrt{24an+25b^2}
\\\iff&\f {144}5 an+36b^2\gs (120a^2+b)^2+24an+25b^2
\\&\qquad\qquad\qquad\ \ \ +2(120a^2+b)\sqrt{24an+25b^2}
\\\iff&\l(\sqrt{24an+25b^2}-5(120a^2+b)\r)^2
\\&\gs 25(2880a^4+48a^2b-b^2+(120a^2+b)^2)=2(60a^2+b)(60a)^2
\end{align*}
and
\begin{align*}
&\sqrt{24an+25b^2}\gs 5(120a^2+b)+60a\sqrt{2(60a^2+b)}
\\\iff& 24an+25b^2\gs 25(120a^2+b+12a\sqrt{2(60a^2+b)})^2
\\\iff&\f n{25}\gs 22a(60a^2+b)+(120a^2+b)\sqrt{2(60a^2+b)}.
\end{align*}
As \eqref{25} holds, by the above the length of the interval
$$J=\l[\f{\sqrt{24an+25b^2}-5b}{12a},\ \f{\sqrt{4an/5+b^2}-b}{2a}\r]$$
is at least $10a$. So, in view of the first paragraph of this proof, there is a positive integer $B\in J$ such that $(n-5bB)/a-5B^2$ is an even integer not in the set $E$ given by \eqref{E}.
Since $B\in J$, we also have
$$0\ls \f{n-5bB}a-5B^2\ls B^2.$$
By Lemma \ref{Lem4.1}, there are $x,y,z\in\Z$ such that
$$\f{n-5bB}a-5B^2=x^2+y^2+z^2+\f{(x+y+z)^2}2=x^2+y^2+z^2+2w^2,$$
where $w=-(x+y+z)/2\in\Z$. As
$$x^2+y^2+z^2+2w^2=\f{n-5bB}a-5B^2\ls B^2,$$
all the four integers
$$x_0=x+B,\ y_0=y+B,\ z_0=z+B,\ w_0=w+B$$
are nonnegative. Note that
\begin{align*}&2w_0(aw_0+b)+x_0(ax_0+b)+y_0(ay_0+b)+z_0(az_0+b)
\\=&\ a(2(w+B)^2+(x+B)^2+(y+B)^2+(z+B)^2)
\\&\ +b(2(w+B)+(x+B)+(y+B)+(z+B))
\\=&\ a(2w^2+x^2+y^2+z^2)+(2aB+b)(2w+x+y+z)+5aB^2+5bB
\\=&(n-5bB-5aB^2)+(2aB+b)\times0+5aB^2+5bB=n.
\end{align*}
This concludes our proof. \qed

\section{Proof of Theorem \ref{Th3}}
\setcounter{lemma}{0}
\setcounter{theorem}{0}
\setcounter{corollary}{0}
\setcounter{remark}{0}
\setcounter{equation}{0}
\setcounter{conjecture}{0}

We need the following lemma  which is essentially \cite[Lemma 3.1]{MS}.

\begin{lemma} \label{LemMS} Let $c,d\in\Z^+$ with $c\eq d\pmod2$, and $c\eq3\pmod9$ or $3\nmid d$.
If $6c>d^2$ and $5c<d^2+2d+6$, then  there are $w,x,y,z\in\N$ such that
\begin{equation}\label{3cd} 3w^2+x^2+y^2+z^2=c
\ \ \t{and}\ \ 3w+x+y+z=d.
\end{equation}
\end{lemma}

Now we give another lemma.

\begin{lemma} \label{Lem3cd} Let $a\in\Z^+$ and $b,c,d\in\Z$ with $b>-a$ and
$ac+bd=n\gs 14a/5+2b$. Then
$$d\gs 3\ \ \t{and}\ \ \f{d^2}6<c<\f{d^2+2d+6}5$$
if and only if
$$\f{\sqrt{5a(n-a+b)+25b^2/4}-a-5b/2}a<d<\f{\sqrt{6an+9b^2}-3b}a.$$
\end{lemma}
\Proof.
Observe that $$3a+\f 52b>3a-\f 52 a=\f a2>0$$ and hence
\begin{align*}&\f{\sqrt{5a(n-a+b)+25b^2/4}-a-5b/2}a\gs2
\\\iff&5a(n-a+b)+\f{25}4b^2\gs \l(3a+\f52b\r)^2
\\\iff&n\gs\f{14}5a+2b.
\end{align*}
Since $n\gs 14a/5+2b$, by the above, if
$$d>\f{\sqrt{5a(n-a+b)+25b^2/4}-a-5b/2}a,$$
then $d >2$ and hence $d\gs3$.

Below we assume $d\gs3$. Then
$$ad+3b>a(d-3)\gs0\ \ \t{and}\ \ ad+a+\f 52b>\f 52(a+b)>0.$$ Observe that
\begin{align*}&\f{d^2}6<c<\f{d^2+2d+6}5
\\\iff& \f{a}6d^2<ac=n-bd<\f a5(d^2+2d+6)
\\\iff&\sqrt{6an+9b^2}>ad+3b\ \t{and}\ \sqrt{5a(n-a+b)+\f{25}4b^2}<ad+a+\f 52b.
\end{align*}
So the desired result follows. \qed

\medskip
\noindent{\tt Proof of Theorem \ref{Th3}}.
By \eqref{n6}, we have
\begin{align*}\f n6&>66a^2(9a-1)+(33a-1)b= 324a^3+(270a^3-66a^2)+(33a-1)b
\\&>a(18a-1)^2+(18a-1)b + (270a^3-66a^2+15ab)
\\&>a(18a-1)^2+(18a-1)b+66a^2(a-1)+15a(a+b)
\\&>a(18a-1)^2+(18a-1)b
\end{align*}
and hence
$$\f{an}6+\f{b^2}4>\l(a(18a-1)+\f b2\r)^2.$$
Let $L$ be the length of the interval
$$K=\l(\f{\sqrt{5a(n-a+b)+25b^2/4}-a-5b/2}a,\,\f{\sqrt{6an+9b^2}-3b}a\r).$$
It is easy to verify that
\begin{align*}&L>18a
\\\iff&\sqrt{6an+9b^2}-\l(a(18a-1)+\f b2\r)>\sqrt{5a(n-a-b)+\f{25}4b^2}
\\\iff&\l(\sqrt{6an+9b^2}-\l(a(18a-1)+\f b2\r)\r)^2>5a(n-a+b)+\f{25}4b^2
\\\iff&\l(\sqrt{6an+9b^2}-3(2a(18a-1)+b)\r)^2>36a^2\times 15(18a^2-2a+b).
\end{align*}
By \eqref{n6}, we have
$$6an+9b^2>\l(3(2a(18a-1)+b)+6a\sqrt{15(2a(9a-1)+b)}\r)^2$$
and hence
$$\l(\sqrt{6an+9b^2}-3(2a(18a-1)+b)\r)^2>36a^2\times 15(18a^2-2a+b).$$
So $L>18a$ by the above.

Now we claim that there are $c,d\in\Z$ with $c\eq d\pmod2$
such that $d\in K$, $ac+bd=n$, and $c\eq3\pmod9$ if $3\mid d$.
To prove this, we distinguish two cases.

{\it Case} 1. $3\nmid a$ or $3\nmid n$.

As $\gcd(a-b,2a)=(a-b,2)$ divides $n$, we may choose the least integer $d_0\in K$
such that $(b-a)d_0\eq n\pmod{2a}$.
If $3\nmid a$, then we choose $d\in\{d_0,d_0+2a\}$ such that $3\nmid d$. If $3\mid a$ and $3\nmid n$, then $d=d_0\not\eq0\pmod3$.
As $L>4a$, we have $d\in K$. Write $n=(b-a)d+2ac_0$ with $c_0\in\Z$.
Then $n=ac+bd$ with $c=2c_0-d\eq d\pmod2$.

{\it Case} 2. $3\mid a$ and $3\mid n$.

As $3\mid a$, we have $3\nmid b$ in this case. Note that
$\gcd(a-b,6a)=(a-b,2)$ divides $n$. As $L>18a$, we may choose the least integer $d_0\in K$ such that
$(b-a)d_0\eq n\pmod{6a}$. Since $3\mid n$ and $3\nmid b-a$, we must have $3\mid d_0$.
Note that $(n-bd_0)/a=-d_0+6q$ for some $q\in\Z$. Choose $r\in\{0,1,2\}$ such that
$r\eq q+1+d_0/3\pmod 3$. Then $d=d_0+6ar\in K$ since $L>18a$. Note that
$$c:=\f{n-bd}a=\f{n-b(d_0+6ar)}a=-d_0+6q-6r\eq -d_0-3q+3r\eq3\pmod9$$
and also $c\eq -d_0\eq d\pmod2$.

By the proved claim, we can write $n$ as $ac+bd$ with $c,d\in \Z$ and $d\in K$
such that $c\eq d\pmod2$, and $c\eq3\pmod 9$ if $3\nmid d$.
Note that
$$\f n6>66a^2(9a-1)+a+(33a-1)(-a)=495a^3+99a^2(a-1)+2a>3a.$$
No matter $b<0$ or $b\gs0$, we always have $n>3a+2b$.
As $d\in K$, by Lemma \ref{Lem3cd} we have $d\gs3$, $6c>d^2$ and $5c<d^2+2d+6$.
Applying Lemma \ref{LemMS}, we see that there are $w,x,y,z\in\N$ satisfying \eqref{3cd}.
Therefore
\begin{align*}n&=ac+bd=a(3w^2+x^2+y^2+z^2)+b(3w+x+y+z)
\\&=3w(aw+b)+x(ax+b)+y(ay+b)+z(az+b).
\end{align*}
This concludes our proof of Theorem \ref{Th3}. \qed

\section{Some remarks and conjectures}
\setcounter{lemma}{0}
\setcounter{theorem}{0}
\setcounter{corollary}{0}
\setcounter{remark}{0}
\setcounter{equation}{0}
\setcounter{conjecture}{0}

It is easy to see that
$$\l\{\f{n(n+1)}2:\ n\in\N\r\}=\{x(2x+1):\ x\in\Z\}.$$
For integers $a>b\gs0$ and $a\eq b\pmod 2$, clearly
$$\l\{\f{x(ax+b)}2:\ x\in\Z\r\}\se\N.$$
The author \cite{S18,S20} investigated those tuples $(a,b,c,d,e,f)\in\Z^6$
with $a>b\gs0$, $c>d\gs0$, $e>f\gs0$ and $a-b,c-d,e-f\in 2\Z^+$
such that each $n\in\N$ can be written as
$$\f{x(ax+b)}2+\f{y(cy+d)}2+\f{z(ez+f)}2$$
with $x,y,z\in\Z$. Some further results along this direction were proved in \cite{WS}.

The following lemma can be easily obtained based on the proof of Lemma \ref{Lem2.1}
given in \cite{MS}.

\begin{lemma} \label{Lem6.1} Let $c,d\in\N$ with $4\nmid c$, $c\eq d\pmod2$  and $4c>d^2$.  Then there are
$s,t,u,v\in\Z$ such that
\begin{equation}\label{12cd} s^2+t^2+u^2+v^2=c\ \t{and}\ s+t+u+v=d.
\end{equation}
\end{lemma}

Using this lemma and some arguments similar to but easier than those in Section 2, we can deduce the following two theorems.

\begin{theorem}\label{Th5.1} Let $a$ and $b$ be positive odd integers with $a>b$ and $\gcd(a,b)=1$.
Then any integer $n>a^3/2+ab$ can be written as
$$\f{w(aw+b)}2+\f{x(ax+b)}2+\f{y(ay+b)}2+\f{z(az+b)}2$$
with $w,x,y,z\in\Z$.
\end{theorem}

\begin{theorem} \label{Th5.2} Let $a,b\in\Z^+$ with $a>b$, $a\not\eq b\pmod2$ and $\gcd(a,b)=1$.

{\rm (i)} If $2\mid a$, $n\in\N$ and $n\gs 2^{1+(-1)^n}a^3+(3+(-1)^n)ab$, then there are $w,x,y,z\in\Z$ such that
$$n=w(aw+b)+x(ax+b)+y(ay+b)+z(az+b).$$

{\rm (ii)} If $2\mid b$, then any integer $n\gs a^3+2ab$ with $4\nmid n$ can be written as
$$w(aw+b)+x(ax+b)+y(ay+b)+z(az+b)$$
with $w,x,y,z\in\Z$.
\end{theorem}

The following conjecture can be viewed as a supplement to Theorem \ref{Th1.3} and Theorem \ref{Th5.2}(ii).

\begin{conjecture} For integers $a>b>0$ with $2\mid b$ and $\gcd(a,b)=1$, if $(a,b)\not=(3,2)$
then there are infinitely many positive multiples of four not representable as
$w(aw+b)+x(ax+b)+y(ay+b)+z(az+b)$
with $w,x,y,z\in\Z$.
\end{conjecture}

In contrast with Corollary \ref{Cor-p57} and Remark \ref{Rem-p57}, we make the following conjecture.

\begin{conjecture} We have
\begin{align*}\{p_5(x)+p_5(y)+p_5(z):\ x,y,z\in\N\}&\supseteq\{n\in\N:\ n>33066\},
\\\{\bar p_5(x)+\bar p_5(y)+\bar p_5(z):\ x,y,z\in\N\}&\supseteq\{n\in\N:\ n>24036\},
\\\{p_6(x)+p_6(y)+p_6(z):\ x,y,z\in\N\}&\supseteq\{n\in\N:\ n>146858\},
\\\{\bar p_6(x)+\bar p_6(y)+\bar p_6(z):\ x,y,z\in\N\}&\supseteq\{n\in\N:\ n>138158\}.
\end{align*}
\end{conjecture}

Motivated by Corollary \ref{Cor5-1}, we pose the following four conjectures.

\begin{conjecture} {\rm (i)} Any integer $n>1$ can be written as
$x(5x+1)+y(5y+1)/2+z(5z+1)/2$ with $x,y,z\in\Z$.

{\rm (ii)} If $n\in\N$ can be written uniquely as $x(5x+1)+y(5y+1)/2+z(5z+1)/2$ with $x,y,z\in\Z$
and $y\ls z$, then it belongs to the set
$$ \{0,\, 2,\, 3,\, 5,\, 7,\, 14,\, 16,\, 19,\, 37,\, 43,\, 58,\, 61,\, 79\}.$$
\end{conjecture}
\begin{remark} For the number of ways to write $n\in\N$ as $x(5x+1)+y(5y+1)/2+z(5z+1)/2$ with $x,y,z\in\Z$
and $y\ls z$, the reader may consult \cite{S24}.
\end{remark}

\begin{conjecture} We have
\begin{align*}\l\{\f{x(5x+1)}2+\f{y(5y+1)}2+\f{z(5z+1)}2:\ x,y,z\in\Z\r\}&=\N\sm\{1,10,19,94\},
\\\l\{x(5x+1)+\f{y(5y+1)}2+\f{3z(5z+1)}2:\ x,y,z\in\Z\r\}&=\N\sm\{1,5,32\},
\\\l\{x(5x+1)+2y(5y+1)+\f{z(5z+1)}2:\ x,y,z\in\Z\r\}&=\N\sm\{1,5,70\},
\\\l\{x(5x+1)+\f{y(5y+1)}2+\f{5z(5z+1)}2:\ x,y,z\in\Z\r\}&=\N\sm\{1,5\}.
\end{align*}
\end{conjecture}

\begin{conjecture} Any integer $n>51$ can be written as
$$\f{w(5w-1)}2+\f{x(5x-1)}2+\f{y(5y+1)}2+\f{z(5z+1)}2$$ with $w,x,y,z\in\N$.
\end{conjecture}

\begin{conjecture} Let
$$N(1)=114862,\ N(-1)=166897,\ N(3)=196987,\ N(-3)=273118.$$
Then, for each $r\in\{\pm1,\pm3\}$, we can write any integer $n>N(r)$ as
$x(5x+r)/2+y(5y+r)/2+z(5z+r)/2$ with $x,y,z\in\N$.
\end{conjecture}

\begin{conjecture} We have
$$\{p_7(x)+2p_7(y)+3p_7(z):\ x,y,z\in\Z\}=\N\sm\{31,\,77\}$$
and
$$\{p_7(x)+2p_7(y)+5p_7(z):\ x,y,z\in\Z\}=\N\sm\{10,\,16\}.$$
\end{conjecture}
\begin{remark} In 2016, the author \cite[Conjecture 5.4(i)]{S16} conjectured that
$$\{p_7(x)+p_7(y)+2p_7(z):\ x,y,z\in\Z\}=\N\sm\{23\}.$$
\end{remark}


\begin{thebibliography}{99}

\bibitem{G} R. K. Guy, {\it Every number is expressible as the sum of how many polygonal numbers?} Amer. Math. Monthly {\bf 101} (1994), 169--172.

\bibitem{KS} D. Krachun and Z.-W. Sun, {\it On sums of four pentagonal numbers with coefficients},
Electron. Res. Arch. {\bf 28} (2020), 559--566.

\bibitem {LW} J. Li and H. Wang, {\it Sums of four polygonal numbers: precise formulas}, arXiv:2405.14710, 2024.

\bibitem{MS} X.-Z. Meng and Z.-W. Sun, {\it Sums of four polygonal numbers with coefficients},
Acta Arith. {\bf 180} (2017), 229--249.

\bibitem{N87}M. B. Nathanson, {\it A short proof of Cauchy's polygonal number theorem}, Proc. Amer. Math. Soc. {\bf 99} (1987), 22--24.
\bibitem {N96} M. B. Nathanson, {\it Additive Number Theory: The Classical Bases}, Grad. Texts in Math., Vol. 164, Springer, New York, 1996.
\bibitem {S15} Z.-W. Sun, {\it On universal sums of polygonal numbers}, Sci. China Math. {\bf 58} (2015), 1367--1396.
\bibitem {S16} Z.-W. Sun, {\it A result similar to Lagrange's theorem}, J. Number Theory {\bf 162} (2016), 190--211.
\bibitem{S18}
Z.-W. Sun, On universal sums $x(ax+b)/2+y(cy +d)/2+z(ez +f)/2$, \textit{Nanjing Univ. J. Math. Biquarterly}, \textbf{35(2)} (2018), 85--199.

\bibitem{S20}
Z.-W. Sun, Universal sums of three quadratic polynomials, \textit{Sci. China Math.}, \textbf{63}  (2020),  no. 3, 501--520.

\bibitem{S24} Z.-W. Sun, Sequence A377224 in OEIS (On-Line Encyclopedia of Integer Sequences),
{\tt http://oeis.org/A377224}.

\bibitem{WS}
H.-L. Wu and Z.-W. Sun, Some universal quadratic sums over the integers, \textit{Electronic Research Archive}, \textbf{27} (2019), 69--87.

\end{thebibliography}
\end{document}